\documentstyle[11pt]{article}

\begin{document}
\newcommand{\ol }{\overline}
\newcommand{\ul }{\underline }
\newcommand{\ra }{\rightarrow }
\newcommand{\lra }{\longrightarrow }
\newcommand{\ga }{\gamma }
\newcommand{\st }{\stackrel }
\newcommand{\scr }{\scriptsize }
\title{\Large\textbf{Some Properties of Finitely Presented Groups with Topological Viewpoints}}
\author{\textbf{Behrooz Mashayekhy \footnote{Correspondence: mashaf@math.um.ac.ir} \ \ and Hanieh Mirebrahimi}\\
Department of Pure Mathematics,\\ Center of Excellence in Analysis on Algebraic Structures,\\
 Ferdowsi University of Mashhad,\\
P. O. Box 1159-91775, Mashhad, Iran}
\date{ }
\maketitle
\begin{abstract}

In this paper, using  some properties of fundamental groups and
covering spaces of connected polyhedra and CW-complexes, we
present topological proof for some famous theorems about finitely
presented groups.
\end{abstract}


\newpage
\section{Introduction}

There are some famous results about subgroups of free groups, free
products and finitely presented groups with complicated group
theoretical proofs. For example, a famous corollary of the
Reidemeister-Schreier rewriting process[3] tells us that every
subgroup of a finitely presented group with finite index is also
finitely presented. In this paper, using some well-known
relationship between covering spaces of connected polyhedra
(simplicial complexes) and their fundamental groups, we intend to
prove some results for
finitely presented groups with a topological approach.\\

\section{ Notation and Preliminaries }

We suppose that the reader is familiar with some well-known notion
such as free groups, free products and presentation in group
theory and simplicial complexes (polyhedra), covering spaces, and
fundamental groups in algebraic topology.\\
\textbf{ Definition 2.1} Let $T$ be a connected simplicial
complex, then $T$ is called a tree if $dim T\leq 1$ and which
contains no circuits. Let $K$ be a connected simplicial complex
with a maximal tree $T$ in $K$. Define a group $G_{K,T}$ with the
following presentation:
$$G_{K,T}=\langle\ (p,q)\in K\ |\ (p,q)\in T,\ (p,q)(q,r)=(p,r)\
if\ \{p,q,r\} \ is\ a\ simplex\ in\ K\ \rangle. $$

The following are some facts in algebraic topology which we need
in the proof of main results.\\
\textbf{Theorem 2.2 ([5]).} Let $K$ be a connected polyhedron with
a base point $p$. Then its fundamental group $\pi_1(K,p)$ is
isomorphic to $G_{K,T}$, where $T$ is a maximal tree in $K$ (note
that we identify the simplicial complex $K$ with its
underlying set the polyhedron $|K|$).\\
\textbf{Corollary 2.3.}  If $K$ is a graph i.e. a connected
$1$-complex, then $\pi_1(K,p)$ is a free group of rank
$|\{(p,q)\in K\backslash
T\ |\ T\ is\ a\ maximal\ tree\ in\ K \}|.$\\
\textbf{ Theorem 2.4 ([5]).} A group $G$ is finitely presented if
and only if there exists a finite connected polyhedron $X$ with $G
\cong \pi_1(X,p).$\\
\textbf{ Theorem 2.5 ([5]).}  For any group $G$, there exists a
CW-complex $K(G)$ with $$\pi_1(K(G))\cong G\ and\ \pi_n(K(G))=1\
for\ all\ n\geq 2.$$ The space $K(G)$ is called
Eilenberg-MacLane space of $G$.\\
\textbf{ Remark 2.6 ([5]).} With respect to the way of
constructing the Eilenberg-MacLane space, generators and relators
of the group $G$ are in one to one corresponding to
$1$-cells and $2$-cells in $K(G)$.\\
\textbf{ Corollary 2.7.} A group $G$ is finitely presented if and
only if the number of $1$-cells and $2$-cells in
it's Eilenberg-MacLane space $K(G)$ is finite.\\
\textbf{ Theorem 2.8 ([1]).} For any group $G$ and its
Eilenberg-MacLane space, $K$ say, we have
$$ H_2(K)\cong M(G), $$
where $M(G)$ is the Schur multiplier of $G$.\\
\textbf{ Lemma 2.9 ([5]).} Let $(\tilde{X},p)$ be a covering space
of $X$, $x_0\in X$, and $Y=p^{-1}(x_0)$ be the fiber over $x_0$.
Then
$|Y|=[\pi_1(X,x_0):p_*(\pi_1(\tilde{X},x_0))]$.\\
\textbf{ Definition 2.10.} A space $X$ is called semilocally
$1$-connected if for every $x\in X$ there exists an open
neighborhood $U$ of $x$ so that every closed path at $x$ in $U$ is
nullhomotopic in $X$.

Note that any CW-complex, particularly any Eilenberg-MacLane
space, is semilocally $1$-connected space.\\
\textbf{ Theorem 2.11 ([5]).} If $X$ is connected, locally path
connected, and semilocally $1$-connected and $G\leq \pi_1(X,x_0)$,
then there exists a constructed covering space of $X$,
$(\tilde{X}_G,p)$ such that
$$ p_*(\pi_1(\tilde{X}_G,\tilde{x_0}))=G. $$
\textbf{ Theorem 2.12 ([5]).} If $X$ is a connected CW-complex and
$\tilde{X}$ is a covering space of $X$, then $\tilde{X}$ is also a
CW-complex with $dim \tilde{X} = dim X$. Moreover, if $X$ has $m$
$k$-cells, and $\tilde{X}$ is $n$-sheeted, then the number of
$k$-cells in $\tilde{X}$ is exactly equal to $mn$.\\
\textbf{Theorem 2.13 ([4]).} For any two groups $G_1$ and $G_2$
with their Eilenberg-MacLane spaces $K_1$ and $K_2$, respectively,
the topological wedge space $K_1\vee K_2$ is an Eilenberg-MacLane
space corresponding to the free product $G_1*G_2$.\\
\textbf{Theorem 2.14 ([4]).} For any two groups $G_1$ and $G_2$
with their Eilenberg-MacLane spaces $K_1$ and $K_2$, respectively,
the topological product space $K_1\times K_2$ is an
Eilenberg-MacLane space corresponding to the direct product
$G_1\times G_2$.

\section{ Main Results}

The following theorem is a consequence of the
Reidemeister-Schreier
rewriting process [3, Prop. 4.2].\\
\textbf{ Theorem 3.1.}\textit{ Every subgroup of a finitely
presented group with finite index is also finitely
presented.}\\
\textbf{ Proof.} Let $G$ be a finitely presented group and $H\leq
G$ with finite index. By Theorem 2.4, there exists a finite
connected polyhedron $X$ with $G\cong \pi_1(X)$. Since $X$ is
connected, locally path connected and semilocally $1$-connected,
there exists a covering space $\tilde{X}_H$ so that
$\pi_1(\tilde{X}_H)\cong H$, by Theorem 2.11. Since $[G:H]\leq
\infty$, $\tilde{X}_H$ is a finite sheeted covering space of $X$
and so by Theorem 2.12, $\tilde{X}_H$ is a finite polyhedron. Now,
by Theorem 2.4, $\pi_1(\tilde{X}_H)\cong H$ is finitely presented.
$\Box$\\
\textbf{ Theorem 3.2.}\textit{ If $G$ is a finitely presented
group, then its Schur multiplier $M(G)$ is finitely
presented. }\\
\textbf{ Proof.} First, note that the Schur multiplier of any
group $G$ is isomorphic to the second homology group of its
corresponding Eilenberg-MacLane space [1], $K$ say. Now using the
fact that the number of $i$-cells, for any $i\in {\bf N}$, in the
Eilenberg-MacLane space of any finitely presented group $G$ is
finite, any homology group of $K$ and in
particular, the Schur multiplier of  $G$ is finitely presented. $\Box$\\
\textbf{ Corollary 3.3.}\textit{ Any covering group
of a finite group is also a finitely presented group.}\\
\textbf{ Proof.} Using the definition of covering group
$\widetilde{G}$ considered as an extension of the Schur multiplier
of $G$ by the group $G$ itself, this note is straightforward
result of two recent theorems. $\Box$\\
\textbf{ Theorem 3.4.}\textit{ The number of
finitely presented groups is countable.}\\
\textbf{ Proof.} First, recall that there exists a bijection
between all finitely generated groups and special $2$-simplicial
complexes [6]. Hence to prove the result, it is sufficient to show
the number of such spaces is countable. Note that each polyhedron
corresponding to a finitely presented group $G$, with a
presentation $G=<x_1,\cdots,x_n\ |\ r_1,\cdots,r_n>$, is obtained
by attaching $r$ $2$-cells to an $n$-rose via some particular
maps.

Suppose $K$ is an $n$-rose lying on the plane ${\bf R}^2$ and
$\{K^n_{\lambda}\}_{\lambda\in \Lambda}$ be the family of all
polyhedra obtained by attaching finitely many $2$-cells to $K$, in
several ways.

Now by Whitney Theorem [7] which states that any $n$-simplicial
complex can be embedded in ${\bf R}^{2n+1}$, we can consider all
the constructed complexes as above in the Euclidean space ${\bf
R}^5$ and then using the axiom of choice and the denseness of
$\textbf{Q}^5$ in $\textbf{R}^5$, we can consider the rational
points $x_{\lambda}\in \textbf{Q}^5$ belonging to one and only one
$K^n_{\lambda}$.

Finally, we conclude that all finitely presented groups with $n$
generators in their presentations are in one to one corresponding
to a subset of rational points in ${\bf R}^5$ and so we are done.
$\Box$\\
\textbf{ Theorem 3.5.}\textit{ The free product of
two finitely presented groups is finitely presented. }\\
\textbf{ Proof.} Suppose that $G_1$ and $G_2$ are finitely
presented groups with Eilenberg-MacLane spaces $K_1$ and $K_2$,
respectively. Using Theorem 2.13, $K_1\vee K_2$ is an
Eilenberg-MacLane space corresponding to $G_1*G_2$. Also, by the
definition, clearly the number of $i$-cells in wedge space of two
spaces having finitely many $i$-cells, is also finite and so by
Theorem 2.7, the result satisfied. $\Box$\\
\textbf{ Theorem 3.6.}\textit{ The product
of two finitely presented groups is finitely presented. }\\
\textbf{ Proof.} By the hypothesis of the previous proof, we only
note Theorem 2.14, and the fact that the number of $i$-cells in
product of two spaces having finitely many $i$-cells, is also
finite. Hence similar to the above proof, we
complete the proof. $\Box$\\
\textbf{ Theorem 3.7.}\textit{ The free amalgamated product of two
finitely presented groups $G_1$ and $G_2$ over a finitely
presented subgroup $H$ is also
finitely presented. }\\
\textbf{ Proof.} First, we consider an Eilenberg-Maclane space
corresponding to the presentation of $H$, $X$ say, and note that
we can extend the algebraic presentation of $H$ to the
presentations for $G_1$ and $G_2$.

Also, by joining some $1$-cells and attaching $2$-cells via the
relations, similar to the method of [5, Theorem 7.45] and [6, Note
6.44], we extend the space $X$ to Eilenberg-Maclane spaces $X_1$
and $X_2$ corresponding to the presentations of $G_1$ and $G_2$,
respectively. Note that the construction is considered so that $X$
is a deformation retract of the space $X_1\cap X_2$.

Now using van-Kampen theorem, the fundamental group $\pi_1(X_1\cup
X_2)$ is the free amalgamated product of two groups
$\pi_1(X_1)\cong G_1$ and $\pi_1(X_2)\cong G_2$ over the subgroup
$\pi_1(Y)\cong \pi_1(X_1\cap X_2)\cong H$ [5].

Hence by uniqueness of the free amalgamated product up to
isomorphism, we conclude that
$$G\cong \pi_1(X_1\cup X_2).$$
On the other hand, by the assumption of being finitely presented
for the groups $H$, $G_1$ and $G_2$ we conclude the spaces $X$,
$X_1$, $X_2$ and so the space $X_1\cup X_2$ have finitely many
cells, which implies the group $\pi_1(X_1\cup X_2)$ to be finitely
presented. $\Box$\\

Finally, by the definition of two new concepts, the Schur
multiplier of a pair and the Schur multiplier of a triple of
groups [2], we conclude the following results. Note that for a
pair of groups $(G,N)$, the natural epimorphism $G\rightarrow G/N$
induces functorially the continuous map $f: K(G)\rightarrow
K(G/N)$. Suppose that $M(f)$ is the mapping cylinder of $f$
containing $K(G)$ as a subspace and is also homotopically
equivalent to the space $K(G/N)$. We take $K(G,N)$ to be the
mapping cone of the cofibration $K(G)\hookrightarrow M(f)$. The
Schur multiplier of the pair $(G,N)$ is considered as the third
homology group of the cofiber space $K(G,N)$.

In addition, we can extend the above notes to a topological
argument for the Schur multiplier of a triple of groups. If we
consider the space $X$ as the cofibration of the natural sequence
$K(G,N)\rightarrow K(G/M,MN/M)$, which is noted by $K(G,M,N)$ [2,
Sec. 6], then the Schur multiplier of the triple $(G,M,N)$ is
defined to be the fourth homology group of the cofiber space
$K(G,M,N)$.\\
\textbf{ Theorem 3.8.}\textit{ The Schur multiplier of
a pair of finitely presented groups is finitely presented.}\\
\textbf{ Proof.} We remark that a mapping cone space obtained from
two spaces having finitely many cells, have also finitely many
cells, which holds the result. $\Box$

Using a similar argument, we establish the following theorem:\\
\textbf{ Theorem 3.9.}\textit{ The Schur multiplier of a triple
of finitely presented groups is finitely presented.}

\end{document}